\pdfoutput=1
\RequirePackage{ifpdf}
\ifpdf 
\documentclass[pdftex]{sigma}
\else
\documentclass{sigma}
\fi

\usepackage{enumerate,epic,mathrsfs}
\usepackage[all,2cell]{xy}

\newtheorem{innercustomthm}{Theorem}

\begin{document}

\newcommand{\dash}[1]{\text{-}\mathrm{#1}}

\newenvironment{customthm}[1]{\renewcommand\theinnercustomthm{#1}\innercustomthm}{\endinnercustomthm}

\newcommand{\arXivNumber}{1409.5717}

\allowdisplaybreaks

\renewcommand{\thefootnote}{$\star$}

\renewcommand{\PaperNumber}{016}

\FirstPageHeading

\ShortArticleName{Extension Fullness of the Categories of Gelfand--Zeitlin and Whittaker Modules}

\ArticleName{Extension Fullness of the Categories\\
of Gelfand--Zeitlin and Whittaker Modules\footnote{This paper is a~contribution to the Special Issue on New Directions
in Lie Theory.
The full collection is available at
\href{http://www.emis.de/journals/SIGMA/LieTheory2014.html}{http://www.emis.de/journals/SIGMA/LieTheory2014.html}}}

\Author{Kevin COULEMBIER~$^\dag$ and Volodymyr MAZORCHUK~$^\ddag$}

\AuthorNameForHeading{K.~Coulembier and V.~Mazorchuk}

\Address{$^\dag$~Department of Mathematical Analysis, Ghent University, Krijgslaan 281, 9000 Gent, Belgium}
\EmailD{\href{coulembier@cage.ugent.be}{coulembier@cage.ugent.be}}
\URLaddressD{\url{http://cage.ugent.be/~coulembier/}}

\Address{$^\ddag$~Department of Mathematics, Uppsala University, Box 480, SE-751 06, Uppsala, Sweden}
\EmailD{\href{mazor@math.uu.se}{mazor@math.uu.se}}
\URLaddressD{\url{http://www2.math.uu.se/~mazor/}}

\ArticleDates{Received September 25, 2014, in f\/inal form February 20, 2015; Published online February 24, 2015}

\Abstract{We prove that the categories of Gelfand--Zeitlin modules of $\mathfrak{g}=\mathfrak{gl}_n$ and Whittaker
modules associated with a~semi-simple complex f\/inite-dimensional algebra $\mathfrak{g}$ are extension full in the
category of all $\mathfrak{g}$-modules.
This is used to estimate and in some cases determine the global dimension of blocks of the categories of
Gelfand--Zeitlin and Whittaker modules.}


\Classification{16E30; 17B10}

\renewcommand{\thefootnote}{\arabic{footnote}}
\setcounter{footnote}{0}

\section{Introduction}

Homological invariants are useful technical tools in modern representation theory.
As classif\/ication of all modules of a~given (Lie) algebra is a~wild problem in almost all non-trivial and interesting
cases (see e.g.~\cite{BKM, Dr,FNP}), the usual ``reasonable'' setup for the study of representations of a~given (Lie)
algebra assumes some f\/ixed subcategory of the category of all modules.
Therefore, the problem to compare homological invariants for a~given category and some of its subcategories is natural
and important.

Given an Abelian category $\mathcal{A}$ and an Abelian subcategory $\mathcal{B}$ of $\mathcal{A}$ such that the natural
inclusion $\mathcal{B}\hookrightarrow \mathcal{A}$ is exact, we say that $\mathcal{B}$ is {\em extension full} in
$\mathcal{A}$ provided that the natural inclusion induces isomorphisms
$\mathrm{Ext}_{\mathcal{B}}^d(M,N)\cong\mathrm{Ext}_{\mathcal{A}}^d(M,N)$ for all $M,N\in \mathcal{B}$ and all $d\geq
0$, see Subsection~\ref{s1.2} for details (we write $X\in\mathcal{C}$ when~$M$ is an object of some category
$\mathcal{C}$).
Extension fullness is a~useful notion which allows one to freely transfer homological information between categories
$\mathcal{A}$ and $\mathcal{B}$.
Recently this concept of extension fullness has also been studied by Herman in~\cite{He}, where it appears under the
name {\em entirely extension closed.}

Motivated by the so-called Alexandru conjecture from~\cite{Fu,Ga} (a part of which asserts extension fullness of certain
categories in Lie theory), in our previous paper~\cite{CoMa} we proved that the category $\mathcal{O}$ associated with
a~semi-simple complex f\/inite-dimensional Lie algebra $\mathfrak{g}$ is extension full in the category of all weight
$\mathfrak{g}$-modules, and that the thick version of $\mathcal{O}$ is extension full in the category of all
$\mathfrak{g}$-modules.
As a~bonus, we determined the global dimension of the thick category $\mathcal{O}$ as well as projective dimensions of
its simple objects.

Although category $\mathcal{O}$ is probably the most studied category of $\mathfrak{g}$-modules, there are several other
natural and well-studied categories which have rather dif\/ferent f\/lavor.
One of them is the category $\mathscr{G}\mathscr{Z}$ of
so-called {\em Gelfand--Zeitlin}\footnote{The surname Zeitlin is spelled {Ce{\u i}tlin} in Russian.
It appeared in dif\/ferent transliterations in Latin script, in particular, as Cetlin, Zetlin, Tzetlin and Tsetlin.
However, it seems that the origin of this surname is the German word ``Zeit'' which justif\/ies our present version.}
modules, introduced in~\cite{DFO1} for the algebra $\mathfrak{sl}_3(\mathbb{C})$, in~\cite{DFO2} for the algebra
$\mathfrak{gl}_n(\mathbb{C})$ and in~\cite{Ma2} for orthogonal Lie algebras.
The category $\mathscr{G}\mathscr{Z}$ can be seen as a~generalization of $\mathcal{O}$ in the sense that it contains
both $\mathcal{O}$ and the thick version of $\mathcal{O}$.
The study of Gelfand--Zeitlin modules attracted considerable attention, see
e.g.~\cite{DFO,FGR,FuOv1,FuOv2,Kh,KM,Ma1,Ma3,MaOv,MaSt,Ov,Ov2,Ra} and references therein.
As far as we know, simple generic Gelfand--Zeitlin modules give the richest known family of simple
$\mathfrak{gl}_n(\mathbb{C})$-modules.
This family depends on $\frac{n(n+1)}{2}$ generic parameters.
The f\/irst main result of this paper is the following statement proved in Section~\ref{s2} (we refer to Sections~\ref{s1}
and~\ref{s2} for more details):

\begin{customthm}{A}
\label{tmain1}
The category $\mathscr{G}\mathscr{Z}$ is extension full in the category of all $\mathfrak{gl}_n$-modules.
\end{customthm}

The added dif\/f\/iculty of the category $\mathscr{G}\mathscr{Z}$ in comparison with thick category $\mathcal{O}$
in~\cite{CoMa} is that~$\mathscr{G}\mathscr{Z}$ is not a~Serre subcategory generated by a~well-known category which has
enough projective objects (contrary to the relation between category~$\mathcal{O}$ and thick category $\mathcal{O}$).
Therefore to prove Theorem~\ref{tmain1} we have to modify and strengthen the abstract results on extension fullness
in~\cite{CoMa}.
Our arguments also heavily use some properties of Gelfand--Zeitlin modules established by Futorny and Ovsienko
in~\cite{FuOv2}.

Another big class of $\mathfrak{g}$-modules, where now $\mathfrak{g}$ is an arbitrary semi-simple Lie algebra with
a~f\/ixed triangular decomposition $\mathfrak{g}=\mathfrak{n}_-\oplus\mathfrak{h}\oplus \mathfrak{n}_+$, is the class of
so-called {\em Whittaker} modules introduced by Kostant in~\cite{Ko}.
Simple Whittaker modules are simple $\mathfrak{g}$-modules on which the algebra $U(\mathfrak{n}_+)$ acts locally
f\/initely, see also~\cite{BM} for a~general Whittaker setup.
These modules were studied in~\cite{KhMa,McD,McD2,MiSo,Se} in the classical setup.
Generalizations of these modules for (inf\/inite-dimensional) Lie algebras and some related algebras attracted a~lot of
attention recently, see~\cite{BCW,BM,BO,Ch,GL,On,WiOn1,WiOn2} and references therein.

We def\/ine the category $\mathscr{W}$ of Whittaker modules of f\/inite length and for this category we prove the next
statement, which is our second main result (we refer to Section~\ref{s3} for more details):

\begin{customthm}{B}
\label{tmain2}
The category $\mathscr{W}$ is extension full in the category of all $\mathfrak{g}$-modules.
\end{customthm}

Using adjunction, the study of $\mathscr{W}$ reduces to the study of locally f\/inite modules over a~certain Noetherian
algebra.
An added dif\/f\/iculty compared to the case of the category $\mathscr{G}\mathscr{Z}$ is that a~module in $\mathscr{W}$ does
not decompose into a~direct sum of f\/inite-dimensional $U(\mathfrak{n}_+)$-modules despite the fact that the action of
$U(\mathfrak{n}_+)$ is locally f\/inite.
Moreover, a~module in $\mathscr{W}$ is not always f\/initely generated over $U(\mathfrak{n}_+)$.
To be able to prove extension fullness, we crucially depend on a~result of Donkin and Dahlberg from the 80's
(see~\cite{Do} and~\cite{Da}) asserting that essential extensions of locally f\/inite modules over solvable f\/inite-dimensional Lie algebras are locally f\/inite.

An advantage of the result in Theorems~\ref{tmain1} and~\ref{tmain2} is that it allows to replace the calculation of
extensions in a~category without projective objects by the calculation of extensions in a~category with enough
projective objects.
In particular, as a~consequence we obtain that the global dimension of the categories $\mathscr{G}\mathscr{Z}$ and
$\mathscr{W}$ equals $\dim\mathfrak{g}$.
As both category $\mathcal{O}$ and its thick version are full subcategories in $\mathscr{G}\mathscr{Z}$ or $\mathscr{W}$
(thick category $\mathcal{O}$ is even a~Serre subcategory, as def\/ined in~\cite[Subsection~10.3.2]{We}), combining
Theorems~\ref{tmain1} and~\ref{tmain2} with the results of~\cite{CoMa} gives, in particular, a~lower bound on the
projective dimension of simple highest weight and Verma modules in the categories $\mathscr{G}\mathscr{Z}$ and $\mathscr{W}$.

The paper is organized as follows: in Section~\ref{s1} we prove some preliminary homological algebra statements,
Section~\ref{s2} deals with the case of Gelfand--Zeitlin modules and Section~\ref{s3} is devoted to the case of
Whittaker modules.

\section{Extension full subcategories}\label{s1}

\subsection{Extensions in Abelian categories}

Let $\mathcal{A}$ be an Abelian category and $M,N\in\mathcal{A}$.
Recall that, for $d\in\mathbb{Z}_{\geq 0}$, the set $\mathrm{Ext}_{\mathcal{A}}^d(M,N)$ of {\em degree~$d$ extensions}
from~$M$ to~$N$ is def\/ined as the set of equivalence classes of exact sequences
\begin{gather*}
\xymatrix{
0\ar[r]&N\ar[r]&X_1\ar[r]&X_2\ar[r]&\cdots\ar[r]&X_d\ar[r]&M\ar[r]&0,
}
\end{gather*}
where all objects and all morphisms are in $\mathcal{A}$, modulo the minimal equivalence relation which contains the
binary relation given by existence of a~commutative diagram
\begin{gather*}
\xymatrix{
0\ar[r]&N\ar[r]\ar@{=}[d]&X_1\ar[r]\ar[d]&X_2\ar[r]\ar[d]&\cdots\ar[r]&X_d\ar[r]\ar[d]&M\ar[r]\ar@{=}[d]&0\\
0\ar[r]&N\ar[r]&Y_1\ar[r]&Y_2\ar[r]&\cdots\ar[r]&Y_d\ar[r]&M\ar[r]&0
}
\end{gather*}
The set $\mathrm{Ext}_{\mathcal{A}}^d(M,N)$ has the natural structure of an Abelian group via the Baer sum.
If $\mathcal{A}$ is $\Bbbk$-linear for some f\/ield $\Bbbk$, then $\mathrm{Ext}_{\mathcal{A}}^d(M,N)$ has the structure of
a~$\Bbbk$-vector space.
We refer to~\cite[Section~3.4]{We} for further information and details.

For $M\in\mathcal{A}$, the {\em projective dimension} $\mathrm{proj.dim}(M)$ of~$M$ is def\/ined as the maximal~$d$ such
that there is $N\in \mathcal{A}$ with $\mathrm{Ext}_{\mathcal{A}}^d(M,N)\neq 0$.
If such maximal~$d$ does not exist, then $\mathrm{proj.dim}(M):=\infty$.
Dually one def\/ines the {\em injective dimension} $\mathrm{inj.dim}(N)$ for $N\in\mathcal{A}$.
The {\em global dimension} $\mathrm{gl.dim}(\mathcal{A})\in\mathbb{Z}_{\geq 0}\cup\{\infty\}$ is the supremum of
projective dimensions taken over all objects in $\mathcal{A}$.
The global dimension coincides with the supremum of injective dimensions taken over all objects in $\mathcal{A}$
(see~\cite[Lemma~5.11.11]{Zi}).

\subsection{Extension full subcategories}\label{s1.2}

Let $\mathcal{A}$ be an Abelian category and $\mathcal{B}$ a~full Abelian subcategory of $\mathcal{A}$ in the sense that
the Abelian structure of $\mathcal{B}$ is inherited from $\mathcal{A}$.
In particular, the natural inclusion functor $\boldsymbol{\iota}:\mathcal{B}\to \mathcal{A}$ is exact.
Then, for every $M,N\in\mathcal{B}$ and every $d\in\mathbb{Z}_{\geq0}$, the functor $\boldsymbol{\iota}$ induces
homomorphisms
\begin{gather*}
\iota^d_{M,N}: \ \mathrm{Ext}_{\mathcal{B}}^d(M,N)\to \mathrm{Ext}_{\mathcal{A}}^d(M,N)
\end{gather*}
of Abelian groups.
In general, these homomorphisms $\iota^d_{M,N}$ are neither injective nor surjective.

We say that $\mathcal{B}$ is {\em extension full} in $\mathcal{A}$ provided that $\iota^d_{M,N}$ are bijective for all
$d\in\mathbb{Z}_{\geq0}$ and for all $M,N\in\mathcal{B}$.
We refer to~\cite[Section~2]{CoMa} for details.

Let $0\to K\to M\to N\to 0$ be a~short exact sequence in $\mathcal{B}$.
Then, for $Q\in \mathcal{B}$, application of $\mathrm{Hom}_{\mathcal{B}}(Q,{}_-)$ and
$\mathrm{Hom}_{\mathcal{A}}(Q,{}_-)$ to this short exact sequence produces the usual long exact sequences in homology
for the categories $\mathcal{B}$ and $\mathcal{A}$, respectively.
Moreover, the homomor\-phisms~$\iota^d_{Q,{}_-}$ give rise to a~homomorphism between these long exact sequences.
A~similar statement is true for $\mathrm{Hom}_{\mathcal{B}}({}_-,Q)$ and $\mathrm{Hom}_{\mathcal{A}}({}_-,Q)$.

\subsection{Checking extension fullness}

In this section we formulate and prove three propositions which will be useful for our study of extension fullness later
in the paper.
The following statement is a~modif\/ication of~\cite[Lemma~4]{CoMa} which also allows for a~somewhat stronger formulation.

\begin{proposition}
\label{prop1}
Let $\mathcal{A}$ be an Abelian category and $\mathcal{B}$ a~Serre subcategory of $\mathcal{A}$.
Assume that $\mathcal{B}$ contains a~full subcategory $\mathcal{B}_0$ with the following properties:
\begin{enumerate}[$($a$)$]\itemsep=0pt
\item
\label{prop1.1}
Every object in $\mathcal{B}$ is a~quotient of an object in $\mathcal{B}_0$.
\item
\label{prop1.2}
The homomorphisms $\iota^d_{M,N}$ are bijective for all $d\in\mathbb{Z}_{\geq0}$, $M\in\mathcal{B}_0$ and
$N\in\mathcal{B}$.
\end{enumerate}
Then $\mathcal{B}$ is {\em extension full} in $\mathcal{A}$.
\end{proposition}

\begin{proof}
Our proof is similar to that of~\cite[Lemma~4]{CoMa}.
We prove the statement by induction on~$d$.
Since $\mathcal{B}$ is assumed to be a~Serre subcategory of $\mathcal{A}$, it is clear that $\iota^0_{Q,N}$ and
$\iota^1_{Q,N}$ are isomorphisms for all $Q,N\in\mathcal{B}$.

To prove the induction step, for $Q\in\mathcal{B}$ consider a~short exact sequence
\begin{gather*}
0\to K\to M\to Q\to 0
\end{gather*}
with $M\in\mathcal{B}_0$ and $K\in\mathcal{B}$, which exists by condition~\eqref{prop1.1}.
Applying $\mathrm{Hom}_{\mathcal{B}}({}_-,N)$ and $\mathrm{Hom}_{\mathcal{A}}({}_-,N)$ gives for each~$d$ the following
commutative diagram with exact rows:
\begin{gather*}
\xymatrix{
\mathrm{Ext}^{d-1}_{\mathcal{B}}(M,N)\ar[r]\ar[d]_{\iota^{d-1}_{M,N}}
&\mathrm{Ext}^{d-1}_{\mathcal{B}}(K,N)\ar[r]\ar[d]_{\iota^{d-1}_{K,N}}
&\mathrm{Ext}^d_{\mathcal{B}}(Q,N)\ar[r]\ar[d]_{\iota^d_{Q,N}}
& \mathrm{Ext}^d_{\mathcal{B}}(M,N)\ar[d]_{\iota^d_{M,N}}\\
\mathrm{Ext}^{d-1}_{\mathcal{A}}(M,N)\ar[r]&\mathrm{Ext}^{d-1}_{\mathcal{A}}(K,N)\ar[r]
&\mathrm{Ext}^d_{\mathcal{A}}(Q,N)\ar[r]& \mathrm{Ext}^d_{\mathcal{A}}(M,N)
}
\end{gather*}
Now, $\iota^{d-1}_{M,N}$ and $\iota^d_{M,N}$ are isomorphisms by condition~\eqref{prop1.2} and from the induction step
we have that $\iota^{d-1}_{K,N}$ is a~monomorphism.
The injective four lemma (by which we mean the statement of~\cite[Lemma~I.3.3(i)]{Mc}) therefore implies that
$\iota^d_{Q,N}$ is a~monomorphism.

It remains to show that $\iota^d_{Q,N}$ is an epimorphism.
For this we consider the following commutative diagram with exact rows:
\begin{gather*}
\xymatrix{
\mathrm{Ext}^{d-1}_{\mathcal{B}}(K,N)\ar[r]\ar[d]_{\iota^{d-1}_{K,N}}
&\mathrm{Ext}^{d}_{\mathcal{B}}(Q,N)\ar[r]\ar[d]_{\iota^{d}_{Q,N}}
&\mathrm{Ext}^d_{\mathcal{B}}(M,N)\ar[r]\ar[d]_{\iota^d_{M,N}}
& \mathrm{Ext}^d_{\mathcal{B}}(K,N)\ar[d]_{\iota^d_{K,N}}\\
\mathrm{Ext}^{d-1}_{\mathcal{A}}(K,N)\ar[r]&\mathrm{Ext}^{d}_{\mathcal{A}}(Q,N)\ar[r]
&\mathrm{Ext}^d_{\mathcal{A}}(M,N)\ar[r]& \mathrm{Ext}^d_{\mathcal{A}}(K,N)
}
\end{gather*}
As $\iota^{d-1}_{K,N}$ is a~bijection by the induction step, $\iota^d_{M,N}$ is a~bijection by condition~\eqref{prop1.2}
and $\iota^d_{K,N}$ is a~monomorphism by the previous paragraph, the surjective four lemma (by which we mean the
statement of~\cite[Lemma~I.3.3(ii)]{Mc}) implies that $\iota^d_{Q,N}$ is an epimorphism.
This completes the proof.
\end{proof}

The proof of the following proposition is dual to that of~\cite[Corollary~5]{CoMa} concerning acyclicness
(see~\cite[Subsection~2.4.3]{We}) of projective objects.
\begin{proposition}
\label{propinj}
Let $\mathcal{A}$ be an Abelian category and $\mathcal{B}$ a~full Abelian subcategory of $\mathcal{A}$ and assume that
they both have enough injective objects.
If every injective object in $\mathcal{B}$ is acyclic for the functor $\mathrm{Hom}_{\mathcal{A}}(K,{}_-)$ for any
$K\in\mathcal{B}$, then $\mathcal{B}$ is extension full in $\mathcal{A}$.
\end{proposition}
\begin{proof}
Consider $K\in\mathcal{B}$ f\/ixed.
We need to prove that the functor $\mathrm{Ext}^d_{\mathcal{A}}(K,{}_-)$, when restricted to the category $\mathcal{B}$,
is isomorphic to the functor $\mathrm{Ext}^d_{\mathcal{B}}(K,{}_-)$.
We have the obvious isomorphism
\begin{gather*}
\mathrm{Hom}_{\mathcal{B}}(K,{}_-)\cong\mathrm{Hom}_{\mathcal{A}}(K,{}_-)\circ \boldsymbol{\iota},
\end{gather*}
of functors from the category $\mathcal{B}$ to the category $\mathbf{Sets}$.

By assumption, the exact functor $\boldsymbol{\iota}$ maps injective objects in $\mathcal{B}$ to injective objects in
$\mathcal{A}$.
Injective objects in $\mathcal{A}$ are acyclic for the functor $\mathrm{Hom}_{\mathcal{A}}(K,{}_-)$, that is for all
such objects~$I$ we have $\mathrm{Ext}_{\mathcal{A}}^{d}(K,I)=0$ for all $d>0$.
The classical Grothendieck spectral sequence, see~\cite[Section~5.8]{We}, therefore implies the statement.
\end{proof}

Let us f\/ix the following notation: for an associative algebra~$A$ over a~f\/ield $\Bbbk$ denote by $A\dash{Mod}$ the
category of all~$A$-modules.
We also denote by $A\dash{mod}$ the full subcategory of $A\dash{Mod}$ consisting of all f\/initely generated modules.
We denote by $A\dash{lfmod}$ the full subcategory of $A\dash{Mod}$ consisting of all modules on which the action of~$A$
is locally f\/inite.
Finally, we denote by $A\dash{fmod}$ the full subcategory of $A\dash{mod}$ consisting of all f\/inite-dimensional modules.

\begin{proposition}
\label{NewProp}
Consider an associative algebra~$A$, a~full Abelian subcategory $\mathcal{A}$ in $A\dash{Mod}$, and a~full Abelian
subcategory $\mathcal{B}$ of $\mathcal{A}$.
Assume that these data satisfy the following conditions:
\begin{enumerate}[$($a$)$]\itemsep=0pt
\item
$\mathcal{B}$ is a~Serre subcategory of $A\dash{Mod}$.
\item
For every surjective morphism $\alpha: M\twoheadrightarrow N$, with $M\in \mathcal{A}$ and $N\in \mathcal{B}$, there is
a~$Q\in \mathcal{B}$ and an injective morphism $\beta: Q\hookrightarrow M$ such that the composition
$\alpha\circ\beta:Q\to N$ is surjective.
\end{enumerate}
Then $\mathcal{B}$ is extension full in $\mathcal{A}$.
\end{proposition}
\begin{proof}
For $d>0$ let
\begin{gather*}
\xymatrix@C=5mm{
0\ar[r]&M\ar[r]&X_1\ar[r]&X_2\ar[r]&\cdots\ar[r]&X_{d-1}\ar[r]&X_d\ar[r]&N\ar[r]&0
}
\end{gather*}
and
\begin{gather*}
\xymatrix@C=5mm{
0\ar[r]&M\ar[r]&Y_1\ar[r]&Y_2\ar[r]&\cdots\ar[r]&Y_{d-1}\ar[r]&Y_d\ar[r]&N\ar[r]&0
}
\end{gather*}
be exact sequences in $\mathcal{B}$.
Assume that there is a~commutative diagram in $\mathcal{A}$ with exact rows as follows:
\begin{gather}
\label{eq31}
\begin{split}
& \xymatrix@C=5mm{
0\ar[r]&M\ar[r]\ar@{=}[d]&X_1\ar[r]\ar[d]&X_2\ar[r]\ar[d]&\cdots\ar[r]&
X_{d-1}\ar[r]\ar[d]&X_d\ar[r]\ar[d]&N\ar[r]\ar@{=}[d]&0\\
0\ar[r]&M\ar[r]&Z_1\ar[r]&Z_2\ar[r]&\dots\ar[r]&Z_{d-1}\ar[r]&Z_d\ar[r]&N\ar[r]&0\\
0\ar[r]&M\ar[r]\ar@{=}[u]&Y_1\ar[r]\ar[u]&
Y_2\ar[r]\ar[u]&\cdots\ar[r]&Y_{d-1}\ar[r]\ar[u]&Y_d\ar[r]\ar[u]&N\ar[r]\ar@{=}[u]&0
}
\end{split}
\end{gather}
For $i=1,2,\dots,d$, let $Q_i$ denote the submodule of $Z_i$ generated by the images of $Y_i$ and $X_i$.
Since $\mathcal{B}$ is a~Serre subcategory of $A\dash{Mod}$, we have that $Q_i$ belongs to $\mathcal{B}$.
Then diagram~\eqref{eq31} restricts to the commutative diagram
\begin{gather*}
\xymatrix@C=5mm{
0\ar[r]&M\ar[r]\ar@{=}[d]&X_1\ar[r]\ar[d]&X_2\ar[r]\ar[d]&\cdots\ar[r]&
X_{d-1}\ar[r]\ar[d]&X_d\ar[r]\ar[d]&N\ar[r]\ar@{=}[d]&0\\
0\ar[r]&M\ar[r]&Q_1\ar[r]&Q_2\ar[r]&\dots\ar[r]&Q_{d-1}\ar[r]&Q_d\ar[r]&N\ar[r]&0\\
0\ar[r]&M\ar[r]\ar@{=}[u]&Y_1\ar[r]\ar[u]&
Y_2\ar[r]\ar[u]&\cdots\ar[r]&Y_{d-1}\ar[r]\ar[u]&Y_d\ar[r]\ar[u]&N\ar[r]\ar@{=}[u]&0
}
\end{gather*}
in which the complex in the second row might be not exact.
By assumption, there is a~submodule in $Z_d$ which surjects onto~$N$ and is in $\mathcal{B}$.
The sum of that submodule with $Q_d$ is also in $\mathcal{B}$ since~$\mathcal{B}$ is a~Serre subcategory of
$A\dash{Mod}$.
We denote this resulting submodule by~$T_d$.
The kernel~$K_d$ of the surjection $T_d\twoheadrightarrow N$ is also in $\mathcal{B}$ as $\mathcal{B}$ is Abelian.
By the same reasoning there is a~submodule $T_{d-1}$ of $Z_{d-1}$ which is in $\mathcal{B}$, contains $Q_{d-1}$, maps to
$T_d$ and surjects onto~$K_d$.
Proceeding inductively, we construct, for each $i=1,2,\dots,d-2$, a~submodule $T_{i}$ of $Z_{i}$ which is in
$\mathcal{B}$, contains $Q_{i}$, maps to $T_{i+1}$ and surjects onto the kernel of the map from $T_{i+1}$ to $T_{i+2}$.
This gives the commutative diagram
\begin{gather*}
\xymatrix@C=5mm{
0\ar[r]&M\ar[r]\ar@{=}[d]&X_1\ar[r]\ar[d]&X_2\ar[r]\ar[d]&\cdots\ar[r]&
X_{d-1}\ar[r]\ar[d]&X_d\ar[r]\ar[d]&N\ar[r]\ar@{=}[d]&0\\
0\ar[r]&M\ar[r]&T_1\ar[r]&T_2\ar[r]&\dots\ar[r]&T_{d-1}\ar[r]&T_d\ar[r]&N\ar[r]&0\\
0\ar[r]&M\ar[r]\ar@{=}[u]&Y_1\ar[r]\ar[u]&
Y_2\ar[r]\ar[u]&\cdots\ar[r]&Y_{d-1}\ar[r]\ar[u]&Y_d\ar[r]\ar[u]&N\ar[r]\ar@{=}[u]&0
}
\end{gather*}
in $\mathcal{B}$ with exact rows.
The above implies that the natural map
\begin{gather}
\label{eq32}
\mathrm{Ext}_{\mathcal{B}}^d(N,M)\to \mathrm{Ext}_{\mathcal{A}}^d(N,M)
\end{gather}
is injective.

The construction above also says that for an arbitrary exact sequence
\begin{gather*}
\xymatrix@C=5mm{
0\ar[r]&M\ar[r]&Z_1\ar[r]&Z_2\ar[r]&\cdots\ar[r]&Z_{d-1}\ar[r]&Z_d\ar[r]&N\ar[r]&0\\
}
\end{gather*}
in $\mathcal{A}$ with $M,N\in \mathcal{B}$ there is a~commutative diagram
\begin{gather*}
\xymatrix@C=5mm{
0\ar[r]&M\ar[r]&Z_1\ar[r]&Z_2\ar[r]&\cdots\ar[r]&Z_{d-1}\ar[r]&Z_d\ar[r]&N\ar[r]&0\\
0\ar[r]&M\ar[r]\ar@{=}[u]&T_1\ar[r]\ar[u]&T_2\ar[r]\ar[u]&\cdots\ar[r]&
T_{d-1}\ar[r]\ar[u]&T_d\ar[r]\ar[u]&N\ar[r]\ar@{=}[u]&0\\
}
\end{gather*}
with exact rows and such that the second row is in $\mathcal{B}$.
This means that the natural map~\eqref{eq32} is surjective and hence bijective, completing the proof.
\end{proof}

As an immediate corollary from Proposition~\ref{NewProp} we obtain:

\begin{corollary}
For an associative algebra~$A$ we have that
\begin{enumerate}[$($i$)$]\itemsep=0pt
\item $A\dash{fmod}$ is extension full in $A\dash{lfmod}$,
\item $A\dash{mod}$ is extension full in $A\dash{Mod}$ provided that~$A$ is Noetherian.
\end{enumerate}
\end{corollary}

\subsection{Adjunction lemma}

The following statement is standard when dealing with categories with enough projective or injective objects.
We failed to f\/ind it in the literature in the generality we need, so we provide a~proof without the use of projective or
injective objects.

\begin{proposition}[adjunction lemma]
\label{adjlem}
Let $\mathcal{A}$ and $\mathcal{B}$ be two Abelian categories and $(\mathrm{F},\mathrm{G})$ an adjoint pair of exact
functors $\mathrm{F}:\mathcal{A}\to \mathcal{B}$ and $\mathrm{G}:\mathcal{B}\to \mathcal{A}$.
Then for every $d\in\mathbb{Z}_{\geq 0}$, $N\in \mathcal{A}$ and $M\in \mathcal{B}$ there are isomorphisms
\begin{gather*}
\mathrm{Ext}^d_{\mathcal{B}}(\mathrm{F}(N),M)\cong \mathrm{Ext}^d_{\mathcal{A}}(N,\mathrm{G}(M))
\end{gather*}
natural in both~$N$ and~$M$.
\end{proposition}

\begin{proof}
Applying $\mathrm{F}$ to an exact sequence
\begin{gather*}
\xymatrix@C=7mm{
0\ar[r]&\mathrm{G}(M)\ar[r]&X_1\ar[r]&X_2\ar[r]&\cdots\ar[r]&X_{d-1}\ar[r]&X_d\ar[r]&N\ar[r]&0
}
\end{gather*}
in $\mathcal{A}$, gives an exact sequence
\begin{gather*}
\xymatrix@C=5mm{
0\ar[r]&\mathrm{F}\mathrm{G}(M)\ar[r]&\mathrm{F}(X_1)\ar[r]&
\mathrm{F}(X_2)\ar[r]&\cdots\ar[r]&\mathrm{F}(X_{d-1})\ar[r]&\mathrm{F}(X_d)\ar[r]&\mathrm{F}(N)\ar[r]&0
}
\end{gather*}
in $\mathcal{B}$.
Denote by $\mathrm{K}$ the kernel of the adjunction natural transformation $\mathrm{F}\mathrm{G}\to
\mathrm{Id}_{\mathcal{B}}$.
Then we have the following commutative diagram
\begin{gather}
\label{eq11}
\begin{split}
& \xymatrix@C=5mm{
0\ar[r]&\mathrm{K}(M)\ar[r]\ar@{^{(}->}[d]&\mathrm{K}(M)\ar[r]\ar@{^{(}->}[d]&
0\ar[d]\ar[r]&\cdots\ar[r]&0\ar[d]\ar[r]&0\ar[d]\ar[r]&0\ar[d]\ar[r]&0\\
0\ar[r]&\mathrm{F}\mathrm{G}(M)\ar[r]\ar@{->>}[d]&\mathrm{F}(X_1)\ar[r]\ar@{->>}[d]&
\mathrm{F}(X_2)\ar[r]\ar@{=}[d]&\cdots\ar[r]&\mathrm{F}(X_{d-1})\ar[r]\ar@{=}[d]&
\mathrm{F}(X_d)\ar[r]\ar@{=}[d]&\mathrm{F}(N)\ar[r]\ar@{=}[d]&0\\
0\ar[r]&M'\ar[r]\ar@{^{(}->}[d]&X'\ar[r]\ar@{^{(}->}[d]&
\mathrm{F}(X_2)\ar[r]\ar@{=}[d]&\cdots\ar[r]&\mathrm{F}(X_{d-1})\ar[r]\ar@{=}[d]&
\mathrm{F}(X_d)\ar[r]\ar@{=}[d]&\mathrm{F}(N)\ar[r]\ar@{=}[d]&0\\
0\ar[r]&M\ar[r]&X''\ar[r]&
\mathrm{F}(X_2)\ar[r]&\cdots\ar[r]&\mathrm{F}(X_{d-1})\ar[r]&\mathrm{F}(X_d)\ar[r]&\mathrm{F}(N)\ar[r]&0
}
\end{split}
\end{gather}
with exact rows.
Here the homomorphism from the f\/irst to the second row is given by the natural embedding $\mathrm{K}\hookrightarrow
\mathrm{F}\mathrm{G}$ and the third row is just the corresponding cokernel with the morphism from the second to the
third row being the canonical projection.
In particular, $M':=\mathrm{F}\mathrm{G}(M)/\mathrm{K}(M)$.
The homomorphism from $M'$ to~$M$ is the natural inclusion (coming from the def\/inition of $\mathrm{K}$) and, f\/inally,
$X''$ is def\/ined as the push-out and the map to $\mathrm{F}(X_2)$ is given by the universal property of push-outs.
Functoriality of the construction yields a~group homomorphism
\begin{gather*}
\Phi: \ \mathrm{Ext}^d_{\mathcal{A}}(N,\mathrm{G}(M))\to \mathrm{Ext}^d_{\mathcal{B}}(\mathrm{F}(N),M).
\end{gather*}

Applying $\mathrm{G}$ to an exact sequence
\begin{gather*}
\xymatrix@C=7mm{
0\ar[r]&M\ar[r]&Y_1\ar[r]&Y_2\ar[r]&\cdots\ar[r]&Y_{d-1}\ar[r]&Y_d\ar[r]&\mathrm{F}(N)\ar[r]&0
}
\end{gather*}
in $\mathcal{B}$, gives an exact sequence
\begin{gather*}
\xymatrix@C=5mm{
0\ar[r]&\mathrm{G}(M)\ar[r]&\mathrm{G}(Y_1)\ar[r]&\mathrm{G}(Y_2)\ar[r]&
\cdots\ar[r]&\mathrm{G}(Y_{d-1})\ar[r]&\mathrm{G}(Y_d)\ar[r]&\mathrm{G}\mathrm{F}(N)\ar[r]&0
}
\end{gather*}
in $\mathcal{A}$.
Denote by $\mathrm{C}$ the cokernel of the adjunction natural transformation $\mathrm{Id}_{\mathcal{A}}\to
\mathrm{G}\mathrm{F}$.
Then we have the following commutative diagram
\begin{gather*}
\xymatrix@C=5mm{
0\ar[r]&\mathrm{G}(M)\ar[r]\ar@{=}[d]&\mathrm{G}(Y_1)\ar[r]\ar@{=}[d]&\mathrm{G}(Y_2)\ar[r]\ar@{=}[d]&
\cdots\ar[r]&\mathrm{G}(Y_{d-1})\ar[r]\ar@{=}[d]&
Y''\ar[r]\ar@{->>}[d]&N\ar[r]\ar@{->>}[d]&0\\
0\ar[r]&\mathrm{G}(M)\ar[r]\ar@{=}[d]&\mathrm{G}(Y_1)\ar[r]\ar@{=}[d]&\mathrm{G}(Y_2)\ar[r]\ar@{=}[d]&
\cdots\ar[r]&\mathrm{G}(Y_{d-1})\ar[r]\ar@{=}[d]&
Y'\ar[r]\ar@{^{(}->}[d]&N'\ar[r]\ar@{^{(}->}[d]&0\\
0\ar[r]&\mathrm{G}(M)\ar[r]\ar[d]&\mathrm{G}(Y_1)\ar[r]\ar[d]&\mathrm{G}(Y_2)\ar[r]\ar[d]&
\cdots\ar[r]&\mathrm{G}(Y_{d-1})\ar[r]\ar[d]&
\mathrm{G}(Y_d)\ar[r]\ar@{->>}[d]&\mathrm{G}\mathrm{F}(N)\ar[r]\ar@{->>}[d]&0\\
0\ar[r]&0\ar[r]&0\ar[r]&0\ar[r]&
\cdots\ar[r]&0\ar[r]&\mathrm{C}(N)\ar[r]&\mathrm{C}(N)\ar[r]&0
}
\end{gather*}
with exact rows.
Here the homomorphism from the third to the last row is given by the natural projection of
$\mathrm{G}\mathrm{F}\twoheadrightarrow \mathrm{C}$ and the second row is just the corresponding kernel with the
morphism from the second to the third row being the canonical injection.
In particular, $\mathrm{C}(N):=\mathrm{G}\mathrm{F}(M)/N'$.
The homomorphism from~$N$ to $N'$ is the natural surjection and, f\/inally, $Y''$ is def\/ined as the pullback and the map
from $\mathrm{G}(Y_{d-1})$ is given by the universal property of pullbacks.
Functoriality of the construction yields a~group homomorphism
\begin{gather*}
\Psi: \ \mathrm{Ext}^d_{\mathcal{B}}(\mathrm{F}(N),M)\to \mathrm{Ext}^d_{\mathcal{A}}(N,\mathrm{G}(M)).
\end{gather*}

Finally, we claim that~$\Phi$ and~$\Psi$ are inverses of each other.
Consider the following commutative diagram:
\begin{gather}
\label{eq14}
\begin{split}
& \xymatrix@C=3mm{
0\ar[r]&\mathrm{G}(M)\ar[r]\ar[d]&X_1\ar[r]\ar[d]&X_2\ar[r]\ar[d]&\cdots\ar[r]
&X_{d-1}\ar[r]\ar[d]&X_d\ar[r]\ar[d]&N\ar[r]\ar[d]&0\\
0\ar[r]&\mathrm{G}\mathrm{F}\mathrm{G}(M)\ar[r]\ar@{->>}[d]&\mathrm{G}\mathrm{F}(X_1)\ar[r]\ar@{->>}[d]&
\mathrm{G}\mathrm{F}(X_2)\ar[r]\ar@{=}[d]&\cdots\ar[r]&\mathrm{G}\mathrm{F}(X_{d-1})\ar[r]\ar@{=}[d]&
\mathrm{G}\mathrm{F}(X_d)\ar[r]\ar@{=}[d]&\mathrm{G}\mathrm{F}(N)\ar[r]\ar@{=}[d]&0\\
0\ar[r]&\mathrm{G}(M)\ar[r]&\mathrm{G}(X'')\ar[r]&
\mathrm{G}\mathrm{F}(X_2)\ar[r]&\cdots\ar[r]&\mathrm{G}\mathrm{F}(X_{d-1})\ar[r]&
\mathrm{G}\mathrm{F}(X_d)\ar[r]&\mathrm{G}\mathrm{F}(N)\ar[r]&0\\
0\ar[r]&\mathrm{G}(M)\ar[r]\ar@{=}[u]&\mathrm{G}(X'')\ar[r]\ar@{=}[u]&
\mathrm{G}\mathrm{F}(X_2)\ar[r]\ar@{=}[u]&\cdots\ar[r]&\mathrm{G}\mathrm{F}(X_{d-1})\ar[r]\ar@{=}[u]&
Y''\ar[r]\ar[u]&N\ar[r]\ar[u]&0
}
\end{split}
\end{gather}
Here the second row is obtained from the f\/irst one by applying $\mathrm{G}\mathrm{F}$ and the homomorphism from the
f\/irst to the second row is given by adjunction $\mathrm{Id}_{\mathcal{A}}\to\mathrm{G}\mathrm{F}$.
The second and the third rows and the homomorphism between them are given by applying the exact functor $\mathrm{G}$ to
the two middle rows of~\eqref{eq11}.
Note that, by construction, application of $\mathrm{G}$ identif\/ies the two last rows of~\eqref{eq11}, that is
$\mathrm{G}(M')\cong \mathrm{G}(M)$ and $\mathrm{G}(X')\cong \mathrm{G}(X'')$ (this follows from surjectivity of the
natural transformation $\mathrm{G}\mathrm{F}\mathrm{G}\to \mathrm{G}$ given by adjunction).
Consequently, from the adjunction identities we have that the composition of the maps in the f\/irst column is the
identity on $\mathrm{G}(M)$.
Finally, the last row and the homomorphism to the third row are given by the def\/inition of~$\Psi$.
In particular, from the construction we have that the image of~$N$ in the second row (coming from the f\/irst row) and in
the third row (coming from the last row) coincide.
Hence diagram~\eqref{eq14} gives rise to a~diagram
\begin{gather*}
\xymatrix@C=3mm{
0\ar[r]&\mathrm{G}(M)\ar[r]\ar@{=}[d]&X_1\ar[r]\ar[d]&X_2\ar[r]\ar[d]&\cdots\ar[r]
&X_{d-1}\ar[r]\ar[d]&X_d\ar[r]\ar[d]&N\ar[r]\ar[d]&0\\
0\ar[r]&\mathrm{G}(M)\ar[r]&\mathrm{G}(X'')\ar[r]&
\mathrm{G}\mathrm{F}(X_2)\ar[r]&\cdots\ar[r]&\mathrm{G}\mathrm{F}(X_{d-1})\ar[r]&
Q'\ar[r]&N'\ar[r]&0\\
0\ar[r]&\mathrm{G}(M)\ar[r]\ar@{=}[u]&\mathrm{G}(X'')\ar[r]\ar@{=}[u]&
\mathrm{G}\mathrm{F}(X_2)\ar[r]\ar@{=}[u]&\cdots\ar[r]&\mathrm{G}\mathrm{F}(X_{d-1})\ar[r]\ar@{=}[u]&
Y''\ar[r]\ar[u]&N\ar[r]\ar[u]&0
}
\end{gather*}
with exact rows, where $N'$ is the image of~$N$ in $\mathrm{G}\mathrm{F}(N)$ and $Q'$ is the full preimage of $N'$.
Pulling back along the epimorphism $N\twoheadrightarrow N'$ gives a~commutative diagram
\begin{gather*}
\xymatrix@C=3mm{
0\ar[r]&\mathrm{G}(M)\ar[r]\ar@{=}[d]&X_1\ar[r]\ar[d]&X_2\ar[r]\ar[d]&\cdots\ar[r]
&X_{d-1}\ar[r]\ar[d]&X_d\ar[r]\ar[d]&N\ar[r]\ar@{=}[d]&0\\
0\ar[r]&\mathrm{G}(M)\ar[r]&\mathrm{G}(X'')\ar[r]&
\mathrm{G}\mathrm{F}(X_2)\ar[r]&\cdots\ar[r]&\mathrm{G}\mathrm{F}(X_{d-1})\ar[r]&
Q\ar[r]&N\ar[r]&0\\
0\ar[r]&\mathrm{G}(M)\ar[r]\ar@{=}[u]&\mathrm{G}(X'')\ar[r]\ar@{=}[u]&
\mathrm{G}\mathrm{F}(X_2)\ar[r]\ar@{=}[u]&\cdots\ar[r]&\mathrm{G}\mathrm{F}(X_{d-1})\ar[r]\ar@{=}[u]&
Y''\ar[r]\ar[u]&N\ar[r]\ar@{=}[u]&0
}
\end{gather*}
with exact rows.
The last diagram shows that the extensions given by the f\/irst and the last rows coincide, which proves that $\Psi\Phi$
is the identity map.

The claim that $\Phi\Psi$ is the identity map is proved similarly.
This completes the proof.
\end{proof}

\section{Gelfand--Zeitlin modules}\label{s2}

\textbf{Notation.} For a~Lie algebra $\mathfrak{a}$ we denote by $U(\mathfrak{a})$ its universal enveloping algebra and
by $Z(\mathfrak{a})$ the center of $U(\mathfrak{a})$.

\subsection[Gelfand--Zeitlin subalgebra of $\mathfrak{gl}_n$]{Gelfand--Zeitlin subalgebra of $\boldsymbol{\mathfrak{gl}_n}$}

For $k\in\mathbb{Z}_{>0}$ denote by $\mathfrak{g}_k$ the Lie algebra $\mathfrak{gl}_k(\mathbb{C})$.
Set $U_k=U(\mathfrak{g}_k)$ and $Z_k=Z(\mathfrak{g}_k)$.
We consider the usual chain
\begin{gather*}
\mathfrak{g}_1\subset \mathfrak{g}_2\subset \mathfrak{g}_3\subset \dots \subset \mathfrak{g}_n
\end{gather*}
of ``left upper corner'' embeddings of Lie algebras as depicted on the following picture:

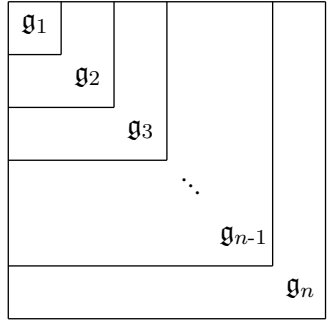
\begin{center}
\begin{picture}(143.00,143.00)
\drawline(20,20)(140,20)
\drawline(20,20)(20,140)
\drawline(140,140)(20,140)
\drawline(140,140)(140,20)
\drawline(20,120)(40,120)
\drawline(40,140)(40,120)
\drawline(20,100)(60,100)
\drawline(60,140)(60,100)
\drawline(20,80)(80,80)
\drawline(80,140)(80,80)
\drawline(20,40)(120,40)
\drawline(120,140)(120,40)
\put(25,130){$\mathfrak{g}_1$}
\put(45,110){$\mathfrak{g}_2$}
\put(65,90){$\mathfrak{g}_3$}
\put(85,70){$\cdot$}
\put(87.5,67.5){$\cdot$}
\put(90,65){$\cdot$}
\put(100,50){$\mathfrak{g}_{n\text{-}1}$}
\put(125,30){$\mathfrak{g}_n$}
\end{picture}
\end{center}

\vspace{-6mm}

This gives rise to the chain
\begin{gather*}
U_1\subset U_2\subset U_3\subset \dots \subset U_n
\end{gather*}
of embeddings of associative algebras.
The subalgebra~$\Gamma$ of $U:=U_n$ generated by all centers~$Z_k$, where $k=1,2,\dots,n$, is called the {\em
Gelfand--Zeitlin} subalgebra.
We set $\mathfrak{g}:=\mathfrak{g}_n$.

The algebra~$\Gamma$ is obviously commutative, moreover,~$\Gamma$ is a~polynomial algebra in $\frac{n(n+1)}{2}$
variables (these can be taken to be generators of $Z_k$ for $k=1,2,\dots,n$), see~\cite{GZ} or~\cite[Chapter~X]{Zh}.
Considered as a~subalgebra of~$U$,~$\Gamma$ is a~{\em Harish-Chandra subalgebra} in the sense of~\cite[Section~1]{DFO},
which means that every f\/initely generated~$\Gamma$-subbimodule of~$U$ is already f\/initely generated both as a~left and
as a~right~$\Gamma$-module, see~\cite[Theorem~24]{DFO}.
Furthermore,~$U$ is free both as a~left and as a~right~$\Gamma$-module, see~\cite{Ov}.
Consequently, the usual induction and coinduction functors
\begin{gather*}
\mathrm{Ind}_{\Gamma}^U:=U\bigotimes_{\Gamma}{}_-
\qquad
\text{and}
\qquad
\mathrm{Coind}_{\Gamma}^U:=\mathrm{Hom}_{\Gamma}(U,{}_-)
\end{gather*}
are exact.

\subsection{Gelfand--Zeitlin modules}

The category $\mathscr{G}\mathscr{Z}$ of {\em Gelfand--Zeitlin}-modules for~$U$ is def\/ined as the full subcategory of
$U\dash{mod}$ (the category of f\/initely generated~$U$-modules), which consists of those $M\in U\dash{mod}$ on which the
action of~$\Gamma$ is locally f\/inite, that is $\dim(\Gamma v)<\infty$ for all $v\in M$.
The category $\mathscr{G}\mathscr{Z}$ is a~Serre subcategory of $U\dash{mod}$ as~$U$ is Noetherian.

Write $\mathrm{Specm}(\Gamma)$ for the set of maximal ideals in~$\Gamma$.
As~$\Gamma$ is a~polynomial algebra, we have the decomposition
\begin{gather*}
\Gamma\dash{fmod}=\bigoplus_{\mathbf{m}\in\mathrm{Specm}(\Gamma)} \Gamma\dash{fmod}_{\mathbf{m}},
\end{gather*}
where $\Gamma\dash{fmod}_{\mathbf{m}}$ denotes the full subcategory of $\Gamma\dash{fmod}$ consisting of all objects
annihilated by some power of $\mathbf{m}$.
From this decomposition, consider the functor
\begin{gather*}
\mathrm{Ind}_{\Gamma,\mathbf{m}}^U: \ \Gamma\dash{fmod}_{\mathbf{m}}\to \mathscr{G}\mathscr{Z}
\end{gather*}
def\/ined as the restriction of $\mathrm{Ind}_{\Gamma}^U$ to $\Gamma\dash{fmod}_{\mathbf{m}}$.

Note that modules in $\mathscr{G}\mathscr{Z}$ are usually inf\/inite-dimensional and hence the usual restriction
$\mathrm{Res}_{\Gamma}^U$ ends up in $\Gamma\dash{lfmod}$ and not in $\Gamma\dash{fmod}$.
However, from~\cite[Corollary~5.3(a)]{FuOv2} it follows that for any $M\in \mathscr{G}\mathscr{Z}$ and
$\mathbf{m}\in\mathrm{Specm}(\Gamma)$ the space
\begin{gather*}
M_{\mathbf{m}}:=\big\{v\in \mathrm{Res}_{\Gamma}^U M: \mathbf{m}^iv=0~\text{for some}~i\big\}
\end{gather*}
is f\/inite-dimensional.
This allows us to def\/ine the functor
\begin{gather*}
\mathrm{Res}_{\Gamma,\mathbf{m}}^U: \ \mathscr{G}\mathscr{Z}\to \Gamma\dash{fmod}_{\mathbf{m}},
\end{gather*}
which sends~$M$ to $M_{\mathbf{m}}$ and is def\/ined on morphisms by restriction.
Then the usual adjunction between induction and restriction implies that $\mathrm{Ind}_{\Gamma,\mathbf{m}}^U$ is left
adjoint to $\mathrm{Res}_{\Gamma,\mathbf{m}}^U$.
The functor $\mathrm{Res}_{\Gamma,\mathbf{m}}^U$ is, in turn, left adjoint to the restriction
\begin{gather*}
\mathrm{Coind}_{\Gamma,\mathbf{m}}^U: \ \Gamma\dash{fmod}_{\mathbf{m}}\to \mathscr{G}\mathscr{Z}
\end{gather*}
of the coinduction functor to $\Gamma\dash{fmod}_{\mathbf{m}}$.

\subsection{The main result}

Our main result in this section is the following:

\begin{theorem}
\label{thm3}
The category $\mathscr{G}\mathscr{Z}$ is extension full in $U\dash{Mod}$.
\end{theorem}

\begin{proof}
To prove Theorem~\ref{thm3}, we would like to apply Proposition~\ref{prop1} for $\mathcal{A}=U\dash{Mod}$,
$\mathcal{B}=\mathscr{G}\mathscr{Z}$ and~$\mathcal{B}_0$ being the full subcategory of $\mathcal{B}$ consisting of
all~$U$-modules isomorphic to $\mathrm{Ind}_{\Gamma}^U N$ for a~f\/inite-dimensional~$\Gamma$-module~$N$.
Taking the assumptions of Proposition~\ref{prop1} into account, Theorem~\ref{thm3} reduces to the following lemma:

\begin{lemma}
Let $Q\in\mathscr{G}\mathscr{Z}$,~$N$ be a~finite-dimensional~$\Gamma$-module and $M=\mathrm{Ind}_{\Gamma}^U N$.
Then $\iota^d_{M,Q}$ is an isomorphism for every $d\in\mathbb{Z}_{\geq 0}$.
\end{lemma}

\begin{proof}
By additivity, we may assume $N\in \Gamma\dash{fmod}_{\mathbf{m}}$ for some $\mathbf{m}\in\mathrm{Specm}(\Gamma)$.
The image of the functor $\mathrm{Ind}_{\Gamma,\mathbf{m}}^U:\Gamma\dash{fmod}_{\mathbf{m}}\to \mathcal{A}$ belongs to
$\mathcal{B}$.
We show that for every $d\in\mathbb{Z}_{\geq 0}$, any $N\in \Gamma\dash{fmod}_{\mathbf{m}}$ and any
$Q\in\mathscr{G}\mathscr{Z}$ we have the isomorphisms
\begin{gather}
\label{eq1}
\mathrm{Ext}^d_{\mathcal{A}}\big(\mathrm{Ind}_{\Gamma,\mathbf{m}}^U N, Q\big)  \cong
\mathrm{Ext}^d_{\Gamma\dash{Mod}}\big(N,\mathrm{Res}_{\Gamma}^U Q\big),
\\
\label{eq2}
\mathrm{Ext}^d_{\Gamma\dash{Mod}}\big(N,\mathrm{Res}_{\Gamma}^U Q\big)  \cong
\mathrm{Ext}^d_{\Gamma\dash{mod}_{\mathbf{m}}}\big(N,\mathrm{Res}_{\Gamma,\mathbf{m}}^U Q\big),
\\
\label{eq3}
\mathrm{Ext}^d_{\Gamma\dash{mod}_{\mathbf{m}}}\big(N,\mathrm{Res}_{\Gamma,\mathbf{m}}^U Q\big)  \cong
\mathrm{Ext}^d_{\mathcal{B}}\big(\mathrm{Ind}_{\Gamma,\mathbf{m}}^U N, Q\big).
\end{gather}
Since, by construction, these three isomorphisms together with the morphism $\iota^d_{M,Q}$ will yield a~commutative
square, we get that $\iota^d_{M,Q}$ is an isomorphism.

Isomorphism~\eqref{eq2} follows from~\cite[Lemma~17]{CoMa}.
Isomorphisms~\eqref{eq1} and~\eqref{eq3} follow from adjunction lemma (Proposition~\ref{adjlem}).
Note that adjunction lemma applies here since~$U$ is free over~$\Gamma$.
\end{proof}

Theorem~\ref{thm3} now follows.
\end{proof}

\subsection[The global dimension of $\mathscr{G}\mathscr{Z}$]{The global dimension of $\boldsymbol{\mathscr{G}\mathscr{Z}}$}

\begin{corollary}
\label{cor5}
We have $\mathrm{gl.dim}(\mathscr{G}\mathscr{Z})=\mathrm{gl.dim}(U\dash{Mod})=\dim\mathfrak{g}$.
\end{corollary}

\begin{proof}
It is well-known, see e.g.~\cite[Corollary~7.7.3]{We}, that the trivial $\mathfrak{g}$-module $\mathbb{C}$ has maximal
possible projective dimension in $U\dash{Mod}$, namely $\dim\mathfrak{g}$.
Since we obviously have $\mathbb{C}\in \mathscr{G}\mathscr{Z}$, the claim follows from Theorem~\ref{thm3}.
\end{proof}

\begin{remark}
The category $\mathscr{G}\mathscr{Z}$ decomposes into a~direct sum of indecomposable {\em blocks},
see~\cite[Theorem~24]{DFO}.
The proof of Corollary~\ref{cor5} implies that the block containing the trivial $\mathfrak{g}$-module has global
dimension $\dim\mathfrak{g}=n^2$.
However, most of the blocks in $\mathscr{G}\mathscr{Z}$, namely all the so-called {\em strongly generic blocks} in the
sense of~\cite[Section~3]{MaOv}, are equivalent to the category of f\/inite dimension modules over $\Gamma_{\mathbf{m}}$,
the completion of~$\Gamma$ with respect to a~maximal ideal $\mathbf{m}$ (this equivalence is induced~by
$\mathrm{Ind}_{\Gamma,\mathbf{m}}^U$), and hence have global dimension $\frac{n(n+1)}{2}$, that is the Krull dimension
of~$\Gamma$.
\end{remark}

\section{Whittaker modules}\label{s3}

\subsection{The category of Whittaker modules}

Let $\mathfrak{g}$ be a~semi-simple f\/inite-dimensional complex Lie algebra with a~f\/ixed triangular decomposition
$\mathfrak{g}=\mathfrak{n}_-\oplus\mathfrak{h}\oplus\mathfrak{n}_+$.
Consider the subalgebra $R=Z(\mathfrak{g})U(\mathfrak{n}_+)$ in $U(\mathfrak{g})$.
Note that~$R$ is not commutative unless $\mathfrak{g}$ is a~direct sum of copies of $\mathfrak{sl}_2$.
Set $U=U(\mathfrak{g})$.

Denote by $\mathscr{W}$ the full subcategory of $U\dash{mod}$ (the category of f\/initely generated~$U$-modules)
consisting of all $\mathfrak{g}$-modules which are locally f\/inite with respect to the action of~$R$ (cf.~\cite[Definition~1.5]{McD}).
Objects in $\mathscr{W}$ will be called {\em Whittaker modules}, which is a~slight modif\/ication of the original notion
from~\cite{Ko}.

\subsection[Simple f\/inite-dimensional~$R$-modules]{Simple f\/inite-dimensional~$\boldsymbol{R}$-modules}

In order to better understand the category $\mathscr{W}$, we start with a~classif\/ication of simple f\/inite-dimensional~$R$-modules.
For this we would need the following fact.

\begin{proposition}
\label{prop21}
The algebra~$R$ is isomorphic to $Z(\mathfrak{g})\otimes_{\mathbb{C}}U(\mathfrak{n}_+)$.
\end{proposition}

\begin{proof}
From the PBW Theorem it follows that the multiplication map
\begin{gather}
\label{eqnnew2}
U(\mathfrak{h})\otimes_{\mathbb{C}}U(\mathfrak{n}_+)\twoheadrightarrow U(\mathfrak{h})U(\mathfrak{n}_+)
\end{gather}
is bijective.
Injectivity of the Harish-Chandra homomorphism $Z(\mathfrak{g})\to U(\mathfrak{h})$ yields that the
$U(\mathfrak{h})$-components of dif\/ferent elements in $Z(\mathfrak{g})$ are dif\/ferent.
Hence~\eqref{eqnnew2} implies that the surjective homomorphism
$Z(\mathfrak{g})\otimes_{\mathbb{C}}U(\mathfrak{n}_+)\twoheadrightarrow R$ given by multiplication is, in fact,
injective, and hence an isomorphism (see also~\cite[Section~3.3]{Ko}).
\end{proof}

Fix a~maximal ideal $\mathbf{m}$ in $Z(\mathfrak{g})$ and a~linear map
$\chi:\mathfrak{n}_+/[\mathfrak{n}_+,\mathfrak{n}_+]\to\mathbb{C}$ and denote by $V_{\mathbf{m},\chi}$ the space
$\mathbb{C}$ endowed with the action of $Z(\mathfrak{g})$ via the projection $Z(\mathfrak{g})/\mathbf{m}\cong
\mathbb{C}$ and with the action of $U(\mathfrak{n}_+)$ via~$\chi$.
The following statement shows that this gives a~complete and irredundant list of pairwise non-isomorphic simple f\/inite-dimensional~$R$-modules.

{\samepage \begin{proposition}\qquad

\begin{enumerate}[$($i$)$]\itemsep=0pt
\item
\label{prop22.1}
For each $\mathbf{m}$ and~$\chi$ as above, $V_{\mathbf{m},\chi}$ is a~simple~$R$-module.
\item
\label{prop22.2}
Each simple finite-dimensional~$R$-module is isomorphic to $V_{\mathbf{m},\chi}$ for some $\mathbf{m}$ and~$\chi$ as
above.
\item
\label{prop22.3}
We have $V_{\mathbf{m},\chi}\cong V_{\mathbf{m}',\chi'}$ if and only if $\mathbf{m}=\mathbf{m}'$ and $\chi=\chi'$.
\end{enumerate}
\end{proposition}

}

\begin{proof}
We have that $V_{\mathbf{m},\chi}$ is a~$Z(\mathfrak{g})\otimes_{\mathbb{C}}U(\mathfrak{n}_+)$-module by construction.
Hence claim~\eqref{prop22.1} follows from Proposition~\ref{prop21}.
Claim~\eqref{prop22.3} is clear by construction.

To prove claim~\eqref{prop22.2}, we note that~$R$ is a~f\/initely generated complex algebra and hence every
simple~$R$-module admits a~central character by Dixmier's version of Schur's lemma, see~\cite[Proposition~2.6.8]{Di}.
Therefore, from Proposition~\ref{prop21} it follows that simple f\/inite-dimensional~$R$-modules are exactly simple f\/inite-dimensional $U(\mathfrak{n}_+)$-modules with the action of $Z(\mathfrak{g})$ given via the projection
$Z(\mathfrak{g})/\mathbf{m}\cong \mathbb{C}$ for some $\mathbf{m}$.
Since $\mathfrak{n}_+$ is nilpotent, all simple f\/inite-dimensional $\mathfrak{n}_+$ have dimension one by Lie's theorem,
see~\cite[Corollary~1.3.13]{Di}.
The claim follows.
\end{proof}

\subsection[The categories $R\dash{fmod}$ and $R\dash{lfmod}$]{The categories $\boldsymbol{R\dash{fmod}}$ and $\boldsymbol{R\dash{lfmod}}$}

For a~maximal ideal $\mathrm{m}$ in $Z(\mathfrak{g})$ and a~linear map
$\chi:\mathfrak{n}_+/[\mathfrak{n}_+,\mathfrak{n}_+]\to\mathbb{C}$ denote by $R\dash{fmod}_{\mathbf{m},\chi}$ the full
subcategory of $R\dash{fmod}$ consisting of all modules for which all simple composition subquotients are isomorphic to
$V_{\mathbf{m},\chi}$.
Def\/ine similarly the subcategory $R\dash{lfmod}_{\mathbf{m},\chi}$ of $R\dash{lfmod}$.

\begin{proposition}
We have decompositions
\begin{gather*}
R\dash{fmod}\cong\bigoplus_{\mathbf{m},\chi}R\dash{fmod}_{\mathbf{m},\chi}
\qquad
\text{and}
\qquad
R\dash{lfmod}\cong\bigoplus_{\mathbf{m},\chi}R\dash{lfmod}_{\mathbf{m},\chi}.
\end{gather*}
\end{proposition}

\begin{proof}
To prove the claim it is suf\/f\/icient to check that
\begin{gather}
\label{eq6}
\mathrm{Ext}_R^1(V_{\mathbf{m},\chi},V_{\mathbf{m}',\chi'})=0
\end{gather}
unless $\mathbf{m}=\mathbf{m}'$ and $\chi=\chi'$.

If $\mathbf{m}\neq\mathbf{m}'$, then~\eqref{eq6} is clear as $Z(\mathfrak{g})$ is central in~$R$.
Assume $\chi\neq \chi'$ and consider some short exact sequence
\begin{gather}
\label{eq7}
0\to V_{\mathbf{m},\chi'}\to M\to V_{\mathbf{m},\chi}\to 0.
\end{gather}
As $\chi\neq\chi'$, there is $a\in \mathfrak{n}_+$ whose action on~$M$ has two dif\/ferent eigenvalues $\chi(a)\neq
\chi'(a)$.
Let~$v$ and~$w$ be the corresponding non-zero eigenvectors in~$M$ (note that each of them is unique up to a~non-zero
scalar).
Then $\mathbb{C}w$ coincides with the image of $V_{\mathbf{m},\chi'}$ in~$M$.
We claim that $\mathbb{C}v$ is an~$R$-submodule (which would mean that the short exact sequence~\eqref{eq7} splits thus
completing the proof).

Indeed, we have $Z(\mathfrak{g})\mathbb{C}v\subset \mathbb{C}v$ as $\mathbb{C}v=\{x\in M: ax=\chi(a)x\}$ and
$Z(\mathfrak{g})$ commutes with~$a$.
Consider some f\/iltration
\begin{gather*}
0=X_0\subset X_1\subset \dots \subset X_k=\mathfrak{n}_+
\end{gather*}
such that $[a,X_i]\subset X_{i-1}$ for all~$i$.
Such f\/iltration exists as $\mathfrak{n}_+$ is nilpotent.
The elements of $X_1$ commute with~$a$ and hence $X_1\mathbb{C}v\subset \mathbb{C}v$.
Let $b\in X_2$ and assume that $bv=\alpha v+\beta w$.
Then, on the one hand,
\begin{gather*}
abv=\chi(a)\big(\alpha v+\beta w\big)+ (\chi'(a)-\chi(a))\beta w= \chi(a)bv+ (\chi'(a)-\chi(a))\beta w.
\end{gather*}
On the other hand, $abv=\chi(a)bv+[a,b]v$.
As $[a,b]\in X_1$, $X_1v\subset \mathbb{C}v$ and~$v$ and~$w$ are linearly independent, we get
$[a,b]v=(\chi'(a)-\chi(a))\beta w=0$, that is $\beta=0$.
Therefore $bv\in \mathbb{C}v$ and thus $X_2v\subset \mathbb{C}v$.
Proceeding inductively, we get $\mathfrak{n}_+v\subset \mathbb{C}v$ and the proof is complete.
\end{proof}

Note that the algebra~$R$ is isomorphic to the enveloping algebra of the nilpotent Lie algebra $\mathfrak{n}_+\oplus
\mathfrak{a}$, where $\mathfrak{a}$ is an Abelian Lie algebra of dimension $\dim(\mathfrak{h})$, see
Proposition~\ref{prop21}.
In particular $R\cong U(\mathfrak{n})$ for some solvable Lie algebra $\mathfrak{n}$.
Therefore we will be able to make use of the following result taken from~\cite[Proposition~1]{Da}, see also~\cite{Do} and~\cite{Fe}.

\begin{lemma}
\label{DaDo}
Let $\mathfrak{n}$ be a~solvable Lie algebra.
Any module $V\in U(\mathfrak{n})\dash{lfmod}$ has an injective hull~$I_V$ in $U(\mathfrak{n})\dash{Mod}$, moreover,
$I_V\in U(\mathfrak{n})\dash{lfmod}$.
\end{lemma}

We note that existence of injective hulls in $U(\mathfrak{n})\dash{Mod}$ is due, in much bigger generality, to
Baer~\cite{Ba}, see also~\cite{ES}.
Lemma~\ref{DaDo} implies the following.

\begin{proposition}
Let $\mathfrak{n}$ be a~finite-dimensional solvable Lie algebra.
Then the category $U(\mathfrak{n})\dash{lfmod}$ has enough injective objects and, moreover,
$U(\mathfrak{n})\dash{lfmod}$ is extension full in $U(\mathfrak{n})$-{\rm Mod}.
\end{proposition}

\begin{proof}
As $U(\mathfrak{n})\dash{lfmod}$ is a~Serre subcategory of $U(\mathfrak{n})\dash{Mod}$, the injective hulls in
$U(\mathfrak{n})\dash{Mod}$ of Lemma~\ref{DaDo} are automatically injective hulls in $U(\mathfrak{n})\dash{lfmod}$.
This proves that $U(\mathfrak{n})\dash{lfmod}$ has enough injective objects.

As $U(\mathfrak{n})\dash{lfmod}$ is a~locally Noetherian Grothendieck category (see~\cite[Appendix~A]{Kr} or~\cite{Ro}),
it follows that each injective object in this category is a~coproduct of indecomposable injective objects and this
decomposition is unique up to isomorphism.
Since $\mathfrak{n}$ is f\/inite-dimensional, $U(\mathfrak{n})$ is Noetherian and hence a~coproduct of injective objects
in $U(\mathfrak{n})\dash{Mod}$ is injective, see~\cite[Proposition~1.2]{Mt}.
This implies that any injective object in $U(\mathfrak{n})\dash{lfmod}$ is also injective when regarded as a~module in
$U(\mathfrak{n})\dash{Mod}$.
The extension fullness thus follows from Proposition~\ref{propinj}.
\end{proof}

As a~consequence, we obtain the following statement.

\begin{corollary}
\label{prop26}
The category $R\dash{lfmod}$ is extension full in $R\dash{Mod}$.
\end{corollary}

\subsection[$R$ versus $U(\mathfrak{g})$]{$\boldsymbol{R}$ versus $\boldsymbol{U(\mathfrak{g})}$}

\begin{proposition}
\label{prop27}
The algebra $U(\mathfrak{g})$ is free both as a~left and as a~right~$R$-module.
\end{proposition}

\begin{proof}
We prove the statement for the right module structure.
The claim for the left module structure then follows by applying the canonical antiautomorphism of $U(\mathfrak{g})$
generated by $x\mapsto -x$ for $x\in \mathfrak{g}$.

Choose a~basis $Y_1,Y_2,\dots,Y_k$ in $\mathfrak{n}_-$, a~basis $H_1,H_2,\dots,H_l$ in $\mathfrak{h}$ and a~basis
$X_1,X_2,\dots,X_k$ in~$\mathfrak{n}_+$.
Then
\begin{gather}
\label{eq33}
\big\{Y_1^{s_1}Y_2^{s_2}\cdots Y_k^{s_k}H_1^{t_1}H_2^{t_2}\cdots H_l^{t_l}X_1^{r_1}X_2^{r_2}\cdots X_k^{r_k}:
s_i,t_i,r_i\in\mathbb{Z}_{\geq0}\big\}
\end{gather}
is a~basis in $U(\mathfrak{g})$ by the PBW theorem.
For an element~$v$ in the above basis, the degree of of~$v$ is the degree of its $\mathfrak{n}_+$-component
$X_1^{r_1}X_2^{r_2}\cdots X_k^{r_k}$, where the grading on $U(\mathfrak{n}_+)$ is def\/ined by giving constants degree $0$
and simple root vectors degree~$1$.
This induces an $\mathbb{Z}_{\geq 0}$-grading on the vector space $U(\mathfrak{g})$.

Choose a~basis $f_1,f_2,\dots,f_m$ of $U(\mathfrak{h})$ as the free module over the algebra of invariants in~$U(\mathfrak{h})$ with respect to the dot-action of the Weyl group~$W$ of $\mathfrak{g}$.
This algebra of invariants is exactly the image of $Z(\mathfrak{g})$ under the Harish-Chandra isomorphism, moreover,
$m=|W|$.
We f\/ix free generators of this algebra of invariants as $p_1,p_2,\dots, p_l$.
Then we choose free generators $z_1,z_2,\dots,z_l$ of $Z(\mathfrak{g})$ as a~polynomial algebra by requiring $z_i-p_i\in
U(\mathfrak{g})\mathfrak{n}_{+}$.
This gives that the set
\begin{gather*}
\underline{\mathbf{B}}:=\big\{Y_1^{s_1}Y_2^{s_2}\cdots Y_k^{s_k} f_hp_1^{t_1}p_2^{t_2}\cdots p_l^{t_l} X_1^{r_1}X_2^{r_2}\cdots
X_k^{r_k}: s_i,t_i,r_i\in\mathbb{Z}_{\geq0},\, h=1,\dots,m\big\}
\end{gather*}
is a~basis of $U(\mathfrak{g})$.

Consider the set
\begin{gather*}
\mathbf{B}:=\big\{Y_1^{s_1}Y_2^{s_2}\cdots Y_k^{s_k}f_j: s_i\in\mathbb{Z}_{\geq0},\, j=1,2,\dots,m\big\}.
\end{gather*}
By Proposition~\ref{prop21}, the set
\begin{gather*}
\mathbf{B}':=\big\{z_1^{q_1}z_2^{q_2}\cdots z_l^{q_l}X_1^{r_1}X_2^{r_2}\cdots X_k^{r_k}: q_i,r_i\in\mathbb{Z}_{\geq0}\big\}
\end{gather*}
is a~basis of~$R$.
Multiplying $Y_1^{s_1}Y_2^{s_2}\cdots Y_k^{s_k}f_j$ on the right with $z_1^{q_1}z_2^{q_2}\cdots
z_l^{q_l}X_1^{r_1}X_2^{r_2}\cdots X_k^{r_k}$ gives
\begin{gather*}
Y_1^{s_1}Y_2^{s_2}\cdots Y_k^{s_k}f_jz_1^{q_1}z_2^{q_2}\cdots z_l^{q_l}X_1^{r_1}X_2^{r_2}\cdots X_k^{r_k}.
\end{gather*}
This can be written as the sum of
\begin{gather}
\label{eq9745}
Y_1^{s_1}Y_2^{s_2}\cdots Y_k^{s_k}f_j p_1^{q_1}p_2^{q_2}\cdots p_l^{q_l} X_1^{r_1}X_2^{r_2}\cdots X_k^{r_k}
\end{gather}
plus a~linear combination of elements~$v$ from the basis~\eqref{eq33} that have strictly higher degree than the degree
of the homogeneous element~\eqref{eq9745}.
Consider the multiplication map
\begin{gather*}
\mathrm{Span}(\mathbf{B})\otimes_{\mathbb{C}} R\to U(\mathfrak{g}).
\end{gather*}
Using the induction on the total degree of the Cartan part, it is easy to check that this map is surjective.
At the same time, we have the induced map
\begin{gather*}
\varphi: \ \mathbf{B}\times \mathbf{B}'\to U(\mathfrak{g}).
\end{gather*}
For $x\in \mathbf{B}$ and $y\in \mathbf{B}'$, by taking the unique element from~\eqref{eq33} of minimal degree which
appears in the expression of $\varphi(x,y)$ with a~non-zero coef\/f\/icient, induces a~bijection from $\mathbf{B}\times
\mathbf{B}'$ to $\underline{\mathbf{B}}$.
This implies that all elements in $\mathbf{B}\mathbf{B}'$ are linearly independent and hence $\mathbf{B}$ is a~basis of
$U(\mathfrak{g})$ as a~free right~$R$-module.
\end{proof}

\begin{proposition}
\label{prop28}
For every finite-dimensional~$R$-module~$V$ the induced $\mathfrak{g}$-module $M(V):=\mathrm{Ind}^{U(\mathfrak{g})}_R V$
has finite length.
\end{proposition}

\begin{proof}
It is enough to prove the claim for $V=V_{\mathbf{m},\chi}$, where $\mathbf{m}$ is a~maximal ideal in $Z(\mathfrak{g})$
and $\chi:\mathfrak{n}_+/[\mathfrak{n}_+,\mathfrak{n}_+]\to\mathbb{C}$.
In this case the statement follows from~\cite[Theorem~2.8]{McD}.
\end{proof}

\begin{corollary}
Every object in $\mathscr{W}$ has finite length.
\end{corollary}

{\sloppy
\begin{proof}
Each $M\in\mathscr{W}$ is generated, as a~$\mathfrak{g}$-module, by some f\/inite-dimensional~$R$-submo\-du\-le~$V$.
By adjunction,~$M$ is thus a~quotient of $M(V)$ and hence the claim follows from Proposi\-tion~\ref{prop28}.
\end{proof}

}

\subsection{The main result}

Our main result in this section is the following:

\begin{theorem}
\label{thm31}
The category $\mathscr{W}$ is extension full in $U\dash{Mod}$.
\end{theorem}

\begin{proof}
We start by proving extension fullness of the category $\widehat{\mathscr{W}}$ which is def\/ined as the full subcategory
of $U\dash{Mod}$ consisting of all modules which are locally~$R$-f\/inite.
The dif\/ference between $\widehat{\mathscr{W}}$ and $\mathscr{W}$ is that we drop the condition of being f\/initely
generated.

We apply Proposition~\ref{prop1} for $\mathcal{A}=U\dash{Mod}$, $\mathcal{B}=\widehat{\mathscr{W}}$ and $\mathcal{B}_0$
being the full subcategory of $\mathcal{B}$ consisting of all~$U$-modules isomorphic to $M(V)$ for some $V\in
R\dash{lfmod}$.

\begin{lemma}
\label{lem32}
Let $Q\in\widehat{\mathscr{W}}$ and $V\in R\dash{lfmod}$.
Then $\iota^d_{M(V),Q}$ is an isomorphism for every $d\in\mathbb{Z}_{\geq 0}$.
\end{lemma}

\begin{proof}
The image of the functor
\begin{gather*}
\mathrm{Ind}_{R}^{U}: \ R\dash{lfmod}\to \mathcal{A}
\end{gather*}
belongs to $\mathcal{B}$ while the image of the functor
\begin{gather*}
\mathrm{Res}_{R}^{U}: \ \mathcal{B}\to R\dash{Mod}
\end{gather*}
belongs to $R\dash{lfmod}$.
Therefore, to prove our lemma, it is enough to show that for every $d\in\mathbb{Z}_{\geq 0}$, any $N\in R\dash{lfmod}$
and any $Q\in\mathscr{W}$ we have the isomorphisms
\begin{gather}
\label{eq71}
\mathrm{Ext}^d_{\mathcal{A}}\big(\mathrm{Ind}_{R}^U N, Q\big)  \cong  \mathrm{Ext}^d_{R\dash{Mod}}\big(N,\mathrm{Res}_{R}^U Q\big),
\\
\label{eq72}
\mathrm{Ext}^d_{R\dash{Mod}}\big(N,\mathrm{Res}_{R}^U Q\big)  \cong  \mathrm{Ext}^d_{R\dash{lfmod}}\big(N,\mathrm{Res}_{R}^U Q\big),
\\
\label{eq73}
\mathrm{Ext}^d_{R\dash{lfmod}}\big(N,\mathrm{Res}_{R}^U Q\big)  \cong  \mathrm{Ext}^d_{\mathcal{B}}\big(\mathrm{Ind}_{R}^U N, Q\big).
\end{gather}
Isomorphisms~\eqref{eq71} and~\eqref{eq73} follow from adjunction lemma (Proposition~\ref{adjlem}) which applies thanks
to Proposition~\ref{prop27}, while isomorphism~\eqref{eq72} is Proposition~\ref{prop26}.
The claim follows.
\end{proof}

Lemma~\ref{lem32} and Proposition~\ref{prop1} thus imply that $\widehat{\mathscr{W}}$ is extension full in
$U\dash{Mod}$.
Hence, to complete the proof of Theorem~\ref{thm31}, it remains to note that, by Proposition~\ref{NewProp}, the category
$\mathscr{W}$ is extension full in $\widehat{\mathscr{W}}$.
\end{proof}

\subsection[The global dimension of $\mathscr{W}$]{The global dimension of $\boldsymbol{\mathscr{W}}$}

\begin{corollary}
\label{cor551}
We have $\mathrm{gl.dim}(\mathscr{W})=\mathrm{gl.dim}(U\dash{Mod})=\dim\mathfrak{g}$.
\end{corollary}

\begin{proof}
Mutatis mutandis Corollary~\ref{cor5}.
\end{proof}

\begin{remark}
The category $\mathscr{W}$ has a~decomposition
\begin{gather*}
\mathscr{W}\cong \bigoplus_{\mathbf{m},\chi}\mathscr{W}_{\mathbf{m},\chi},
\end{gather*}
where $\mathscr{W}_{\mathbf{m},\chi}$ is the full subcategory consisting of all modules which restrict to
$R\dash{lfmod}_{\mathbf{m},\chi}$, see~\cite[Theorem~9]{BM}.
Corollary~\ref{cor551} says that one of these blocks, namely the one corresponding to the trivial central character and
trivial~$\chi$ has global dimension $\dim\mathfrak{g}$.
This particular block contains many simple objects.
Most of the blocks contain only one simple object and are expected to have smaller global dimension.
Note also that thick category $\mathcal{O}$ is a~Serre subcategory of $\mathscr{W}$.
\end{remark}

\subsection*{Acknowledgements}
KC is a~Postdoctoral Fellow of the Research Foundation~-- Flanders (FWO).
VM is partially supported by the Swedish Research Council.
A~substantial part of this work was done during the visit of the authors to the CRM Thematic Semester ``New Directions
in Lie Theory'' in Montreal.
We thank CRM for hospitality and partial support.
We thank the referees for helpful comments.

\pdfbookmark[1]{References}{ref}
\LastPageEnding

\end{document}